\documentclass[11pt]{article}
\usepackage[a4paper,left=17mm,top=20mm,right=17mm,bottom=25mm]{geometry}

\usepackage{authblk}
\usepackage[numbers, sort&compress]{natbib}
\usepackage[margin=30pt,font=small,labelfont=bf]{caption}
\usepackage{amsmath,amsopn,amsthm,amssymb}

\usepackage{graphicx}
\usepackage{mathptmx} 
\usepackage[mathscr]{euscript}
\usepackage{natbib}
\usepackage{mathtools}
\usepackage{enumitem}
\usepackage{xspace}
\usepackage{scalerel}
\usepackage{bold-extra}

\usepackage[
 	colorlinks=true,
	urlcolor=black,
	linkcolor=blue
]{hyperref}

\usepackage{algorithm}
\usepackage[noend]{algpseudocode}
\algrenewcommand\algorithmicrequire{\textbf{Input:}}
\algrenewcommand\algorithmicensure{\textbf{Output:}}

\newtheorem{theorem}{Theorem}[section]
\newtheorem{lemma}[theorem]{Lemma} 
\newtheorem{corollary}[theorem]{Corollary}
\newtheorem{proposition}[theorem]{Proposition}

\newtheorem{definition}[theorem]{Definition}
\newtheorem{observation}[theorem]{Observation}

\newcommand{\Hasse}[1][]{\mathscr{H}\ifthenelse{\equal{#1}{}}{}{(#1)}}
\newcommand{\hasse}[1][]{\mathscr{H}\ifthenelse{\equal{#1}{}}{}{(#1)}}
\DeclareMathOperator{\CC}{\mathtt{C}}
\newcommand{\LCA}{\ensuremath{\operatorname{LCA}}}
\newcommand{\lca}{\ensuremath{\operatorname{lca}}}

\DeclareMathOperator{\indeg}{indeg}
\DeclareMathOperator{\outdeg}{outdeg}

\newcommand{\One}{\mathscr{I}^{\scaleto{1}{4pt}}}
\newcommand{\OneCaption}{\mathscr{I}^{\scaleto{1}{5pt}}}

\newcommand{\rel}{\textsc{Rel}\xspace}

\newcommand{\Olca}{$\One$-$\lca$\xspace}

\newcommand{\Olcarel}{$\One$-$\lca$-$\rel$\xspace}

\newcommand{\OlcaTMP}{$\mathscr{I}$-$\lca$\xspace}

\newcommand{\OlcarelTMP}{$\mathscr{I}$-$\lca$-$\rel$\xspace}

\newcommand{\REV}[1]{\begingroup\color{black}#1\endgroup}
\providecommand{\keywords}[1]{\textbf{\textit{Keywords: }} #1}

\title{Characterizing and Transforming DAGs within the $\mathscr{I}$-\REV{lca} Framework}

\author[ ]{Marc Hellmuth}
\author[ ]{Anna Lindeberg}

\affil[ ]{Department of Mathematics, Faculty of Science,
  Stockholm University, SE-10691 Stockholm, Sweden }

\date{\ }

\begin{document}
\sloppy

\maketitle

\abstract{ 
	We explore the connections between clusters and least common ancestors (LCAs) in directed
    acyclic graphs (DAGs), focusing on \REV{the interplay between so-called $\mathscr{I}$-lca-relevant DAGs and DAGs with the $\mathscr{I}$-lca-property.
    Here, $\mathscr{I}$ denotes a set of integers.}
	 In $\mathscr{I}$-lca-relevant DAGs, each vertex
    is the unique LCA for some subset $A$ of leaves of size $|A|\in \mathscr{I}$, \REV{whereas in a DAG with the $\mathscr{I}$-lca-property there exists a unique LCA for every subset $A$ of leaves satisfying $|A|\in \mathscr{I}$.
	 We elaborate on the difference between these two properties}
    and establish their close relationship to pre-$\mathscr{I}$-ary and $\mathscr{I}$-ary set systems. \REV{This, in turn, 
    generalizes results established for (pre-) binary and $k$-ary set systems.}
    Moreover, we build upon recently established results that use a simple operator $\ominus$, enabling the transformation of
	arbitrary DAGs into $\mathscr{I}$-lca-relevant DAGs. This process reduces
    unnecessary complexity while preserving  key structural properties of the original DAG.
	The set $\mathfrak{C}_G$ consists of all clusters in a DAG $G$, where clusters correspond to the
    descendant leaves of vertices. While in some cases $\mathfrak{C}_H = \mathfrak{C}_G$ when
    transforming $G$ into an $\mathscr{I}$-lca-relevant DAG $H$, it often happens that certain clusters in $\mathfrak{C}_G$ do not
    appear as clusters in $H$. To understand this phenomenon in detail, we characterize the subset
    of clusters in $\mathfrak{C}_G$ that remain in $H$ for DAGs $G$ with the
    $\mathscr{I}$-lca-property. Furthermore, we show that the set $W$ of vertices required to
    transform $G$ into $H = G \ominus W$ is uniquely determined for such DAGs. This, in turn, allows
    us to show that the \REV{``shortcut-free'' version of the} transformed DAG $H$ is always a tree or a galled-tree whenever
    $\mathfrak{C}_G$ represents the clustering system of a tree or galled-tree and $G$ has the
    $\mathscr{I}$-lca-property. In the latter case $\mathfrak{C}_H = \mathfrak{C}_G$ always holds.
}

\medskip\noindent

\section{Introduction}

Rooted networks and, more generally, directed acyclic graphs (DAGs), play a crucial role in
mathematical phylogenetics for modeling complex evolutionary relationships that traditional rooted
trees cannot fully capture \cite{Huber2022-gu,HRS:10,Huson:11}. However, DAGs and networks
inferred from genomic data are often highly complex and tangled \cite{DM:09}. Consequently, several methods have been developed to
simplify DAGs or networks in various ways while preserving their most significant features to
highlight a ``trend of evolution'' \cite{Heiss2024,FRANCIS2021107215,HUBER201630,HL:24,BORDEWICH2016114}. 
\REV{A basic example of such a simplification is the ``suppression'' of vertices with in- and out-degree 
one to transform a DAG into a so-called phylogenetic DAG. These superfluous vertices 
are not supported by observable data -- typically the leaves of the DAG, i.e., vertices without \REV{children}.
In a DAG $G$ with leaf set $L(G)$, a least common ancestor
(LCA) of a subset $A \subseteq L(G)$ is a vertex $v$ that is an ancestor of all $x \in A$ and has no
descendant that also satisfies this property. It is evident that vertices with in- and out-degree one
cannot not serve as LCAs for any subset of leaves and are therefore considered redundant in certain contexts.}

\REV{This work connects several fundamental aspects of DAG-based modeling: the simplification of DAGs under LCA-based conditions,
the analysis of clustering systems derived from DAGs, and the relationship between these concepts and well-known graph classes
such as regular DAGs, phylogenetic trees and galled trees. Our approach is rooted in the framework developed in \cite{HL:24}, where the notion of
$\mathscr{I}$-$\lca$-relevance was introduced.}
A DAG $G$ is $\mathscr{I}$-$\lca$-relevant if every vertex is the unique LCA of some set of leaves
$A$ with size $|A| \in \mathscr{I}$ \REV{where $\mathscr{I}$ denotes a set of integers. 
In particular, every vertex in a DAG $G$ is the unique LCA of some set of leaves if and only if $G$
is $\mathscr{I}$-$\lca$-relevant for $\mathscr{I} = \{1,2,\dots,|L(G)|\}$.}

LCA vertices represent ancestral relationships with
clear phylogenetic signals and direct relevance to the observed ancestral
relationships in the data. 
\REV{Restricting attention to specific cardinalities, i.e., to specific sets $\mathscr{I}$, is motivated by biological applications
as, for example, orthology detection. Two genes, represented as leaves $x$ and $y$, are considered orthologs
if and only if their unique LCA $\lca(x,y)$ is labeled as a ``speciation'' event \cite{Jensen2001-lr,Fitch2000,Hellmuth2013}.
In this context, non-leaf vertices that are not LCAs of any pair of genes do not convey information about orthologous
relationships and can be omitted. Therefore, simplifying a DAG to a $\{1,2\}$-$\lca$-relevant structure preserves
all orthology-relevant information while also simplifying the DAG to make it easier to interpret.
}

\REV{
On the other hand, we are also interested in the set of clusters $\mathfrak{C}_G$ of a given DAG $G$,
where each cluster corresponds to the set of descendant leaves of a vertex.
Characterizing the structural properties of DAGs in terms of their cluster systems
is a challenging task, but it has become an important and active area of research in its own right; see e.g.\ \cite{Baroni:05,sem-ste-03a,Huson:08,Huson:11,ALRR:14,Barthelemy:08,Brucker:09,Gambette:12}
and \cite{Hellmuth2023} for a general overview. As shown in \cite[Cor.~4.9]{HL:24} (and restated in Proposition~\ref{prop:Onerel=>regular}), $\mathscr{I}$-$\lca$-relevant DAGs from which all ``shortcut edges'' have been removed are \emph{regular}, that is, 
they are uniquely determined by their set of clusters; cf.\ Definition~\ref{def:regular}.
Hence, a further aim of this paper is to explore the connection between other structural properties of DAGs related to LCAs and their associated set of clusters. 
}

As shown in \cite{HL:24}, any DAG $G$ can be transformed into an $\mathscr{I}$-$\lca$-relevant DAG
by \REV{``replacing'' certain internal vertices using a simple operator $\ominus$. This operator acts on a subset $W$
of non-leaf vertices -- specifically those that do not serve as the unique LCA for any subset of leaves of
size in $\mathscr{I}$ -- and produces a simplified DAG $G \ominus W$ that retains the essential structure of $G$.
However, such sets $W$ are not uniquely determined in general. Even so, for any such set $W$, the resulting
set of clusters of $G \ominus W$, denoted by $\mathfrak{C}_{G \ominus W}$, is always a subset of the
the set of original cluster system $\mathfrak{C}_G$ of $G$.
Thus, different choices of $W$ may yield different cluster systems, all contained within $\mathfrak{C}_G$.
This raises several natural questions: Which clusters are lost in $\mathfrak{C}_G \setminus \mathfrak{C}_{G \ominus W}$?
Which subsets $W$ are of minimum size for achieving $\mathscr{I}$-$\lca$-relevance? Since answering these questions is difficult
in the general case, we focus on DAGs that already satisfy the $\mathscr{I}$-$\lca$-property, i.e.,
DAGs where every subset $A \subseteq L(G)$ with $|A| \in \mathscr{I}$ has a unique LCA.}

As we shall see, in this case, the set $W$ ensuring that $G \ominus W$ is
$\mathscr{I}$-$\lca$-relevant is uniquely determined and there is a close connection between the
remaining clusters in $\mathfrak{C}_{G \ominus W}$ and so-called $\mathscr{I}$-ary set systems.
This, in particular, generalizes results established for binary or $k$-ary clustering systems, as
established in \cite{Barthelemy:08, SCHS:24}. 
\REV{Finally, we show that if a DAG $G$ with the $\mathscr{I}$-$\lca$-property has a cluster system
$\mathfrak{C}_G$ corresponding to that of a rooted tree or a galled tree, then the simplification
$G \ominus W$ preserves all clusters, i.e., $\mathfrak{C}_{G \ominus W} = \mathfrak{C}_G$.
Moreover, in this case, $G \ominus W$ from which all so-called shortcuts have been removed
coincides with the unique rooted tree or galled tree $T$ such that $\mathfrak{C}_T = \mathfrak{C}_G$.}


\section{Basics}
\label{sec:basics}

\REV{To help the reader navigate more easily through the definitions used in this paper, we provide Table~\ref{tab:sum-def-C} summarizing the main notation and concepts.}

\paragraph{Sets and set systems.}
All sets considered here are assumed to be finite. A \emph{set system $\mathfrak{C}$ (on $X$)} is a collection of subsets of $X$.
For  $\mathscr{I}\subseteq \{1,\dots, |X|\}$,  we write $X(\mathscr{I})$ for the set system comprising all
subsets $A$ of $X$ with $|A|\in \mathscr{I}$ elements. We put $2^X \coloneqq  X(\{1,\dots, |X|\})$ to denote 
\REV{the set of all non-empty subsets of $X$.} 
A set system $\mathfrak{C}$ on $X$ is \emph{grounded} if $\{x\}\in \mathfrak{C} $ for all $x\in X$ and
$\emptyset\notin \mathfrak{C}$, while $\mathfrak{C}$ is a \emph{clustering system} if it is grounded
and satisfies  $X\in \mathfrak{C}$. 
Two sets $M$ and $M'$ \emph{overlap} if $M \cap M' \notin \{\emptyset, M, M'\}$.

\REV{The following types of set systems will play a central role in this contribution.}
\begin{definition}\label{def:pre-I-ary}
A set system $\mathfrak{C}$ on $X$ is \emph{pre-$\mathscr{I}$-ary} if, 
for all $A\in X(\mathscr{I})$, there is a unique inclusion-minimal element $C\in \mathfrak{C}$ such that $A\subseteq  C$. 
A pre-$\mathscr{I}$-ary set system  $\mathfrak{C}$ on $X$ is called \emph{$\mathscr{I}$-ary} if all  $C\in \mathfrak{C}$ satisfy
property
\begin{description}
	\item[\textnormal{($\mathscr{I}$-$C$)}] there is some $A\in X(\mathscr{I})$
														such that $C$ is the unique inclusion-minimal element in $\mathfrak{C}$ with $A\subseteq  C$.
\end{description}
\end{definition}
A pre-$\{1,2\}$-ary, resp., $\{1,2\}$-ary set system is known as  \emph{pre-binary}, resp., 
\emph{binary} set system  \cite{Barthelemy:08}. Moreover, pre-$\{1,\dots,k\}$-ary and $\{1,\dots,k\}$-ary set systems
have been studied in \cite{SCHS:24}. \REV{Definition~\ref{def:pre-I-ary} is illustrated in Figure~\ref{fig:exmpl-2lca}; see the figure caption for further details.}

\begin{center}
\begin{table}[t]\footnotesize
  \setlength{\tabcolsep}{8pt} 
  \renewcommand{\arraystretch}{1.4} 
\REV{
  \begin{tabular}{p{3.5cm}|p{9.4cm}l} \hline
    \multicolumn{2}{l}{\em A DAG $G$ is/satisfies}  & Ref. \\  \hline
    \emph{path-cluster-comparability (PCC)} &for all $u,v\in V(G)$,
    $u$ and $v$ are $\preceq_G$-comparable if and only if
    $\CC_G(u)\subseteq \CC_G(v)$ or $\CC_G(v)\subseteq \CC_G(u)$.
    & Def.~\ref{def:PCC}\\
    \emph{regular} & there is a prescribed isomorphism between $G$ and the Hasse diagram $\Hasse[\mathfrak{C}_G]$ of the clusters in $G$ & Def.~\ref{def:regular} \\
   \emph{\OlcaTMP-relevant} (\OlcarelTMP)	&  all vertices $v$ in $G$ are \Olca vertices, i.e., $v = \lca_G(A)$ for some $A\in X(\mathscr{I})$
	& Def.~\ref{def:LCAlcarel} \\
    \emph{\OlcaTMP-property} & $\lca_G(A)$ is well-defined for all $A\in X(\mathscr{I})$.    & Def.~\ref{def:Olca-prop} \\
    \hline 
   \end{tabular} \medskip
   
   \begin{tabular}{p{3.5cm}|p{9.4cm}l}	\hline
        \multicolumn{2}{l}{\em A set system $\mathfrak{C}$ on $X$ is/satisfies}  & Ref. \\  \hline
    \emph{pre-$\mathscr{I}$-ary} & for all $A\in X(\mathscr{I})$, there is a unique inclusion-minimal element $C\in \mathfrak{C}$ such that $A\subseteq  C$ & 
     Def.~\ref{def:pre-I-ary}\\
     \emph{($\mathscr{I}$-$C$)} & there is some $A\in X(\mathscr{I})$
											such that $C$ is the unique inclusion-minimal element in $\mathfrak{C}$ with $A\subseteq  C$. & Def.~\ref{def:pre-I-ary}\\
    
   \emph{$\mathscr{I}$-ary}  	& pre-$\mathscr{I}$-ary and $\mathscr{I}$-$C$. 
	& Def.~\ref{def:pre-I-ary} \\
	\emph{grounded} & $\{x\}\in \mathfrak{C} $ for all $x\in X$ and $\emptyset\notin \mathfrak{C}$ & Sec.~\ref{sec:basics} \\
	\emph{clustering system} & grounded and $X\in \mathfrak{C}$ & Sec.~\ref{sec:basics} \\
	\hline
	  \end{tabular} }
  \caption{\REV{Summary of main definitions used in this paper.}}
  \label{tab:sum-def-C}
\end{table}
\end{center}

\paragraph{Directed Graphs.}
We consider \emph{directed graphs} $G=(V,E)$ with nonempty vertex set $V(G)\coloneqq V$ 
and edge set $E(G)\coloneqq E \subseteq\left(V\times V\right)\setminus\{(v,v) \mid v\in
V\}$.  For directed graphs $G=(V_G,E_G)$ and $H=(V_H, E_H)$, an
\emph{isomorphism between $G$ and $H$} is a bijective map $\varphi\coloneqq V_G\to V_H$ such that 
$(u,v)\in E_G$ if and only if $(\varphi(u),\varphi(v))\in E_H$. 
If such a map exist, then $G$ and $H$ are \emph{isomorphic}, in symbols $G\simeq H$. 

A directed graph $G$ is \emph{phylogenetic} if it does not contain a vertex $v$ such that $\outdeg_G(v)\coloneqq\left|\left\{u\in V \colon (v,u)\in
E\right\}\right|=1$ and
$\indeg_G(v)\coloneqq\left|\left\{u\in V \colon (u,v)\in
E\right\}\right|\leq 1$. i.e., 
a vertex with only one outgoing edge and at most one incoming edge.

An \emph{\REV{(undirected)} $v_1v_n$-path} $P = (V,E)$ has an ordered vertex set $V = \{v_1,v_2,\ldots,v_n\}$ and the edges in $E$
are precisely of one of the form $(v_i,v_{i+1})$ or $(v_{i+1},v_i)$, $i=1,2,\ldots, n-1$. If all edges in $P$ are precisely of the form $(v_i,v_{i+1})$ for each $i=1,2,\ldots, n-1$, then $P$
is called \emph{directed}.
A \emph{directed acyclic graph (DAG)} is a directed graph that does not contain directed cycles, \REV{that is, 
it contains no directed $v_1v_n$-path and the edge $(v_n,v_1)$.}
A directed graph $G$ is \emph{connected} if there exists an $xy$-path between every pair of vertices $x$ and $y$ 
and \emph{biconnected} if the removal of any vertex and its incident edges from $G$ yields a connected graph.

We can associate a partial order $\preceq_G$  on the vertex set $V(G)$ of a DAG $G$, defined by
$v\preceq_G w$ if and only if there is a directed $wv$-path. In
this case, we say that $w$ is an \emph{ancestor} of $v$ and $v$ is a \emph{descendant} of $w$. If $v\preceq_G w$
and $v\neq w$, we write $v\prec_G w$. If $(u,v) \in E(G)$, then $u$ is a \emph{parent} of $v$ and 
$v$ a \emph{child} of $u$.
\REV{A vertex $x$ in $G$ without children, that is, a $\preceq_G$-minimal vertex of $G$, 
is called a \emph{leaf} of $G$.} We denote by
$L(G)\subseteq V(G)$ the \REV{set of leaves of $G$}.
If $L(G)=X$, then $G$ is a \emph{DAG on $X$}. 
A vertex $v\in V(G)$ of a DAG $G$ that is $\preceq_G$-maximal \REV{ or, equivalently, that has no parents,} is called a \emph{root}, and the set of roots of $G$ is denoted by $R(G)$. 
Note that $L(G)\neq \emptyset$ and $R(G)\neq \emptyset$ for all DAGs $G$ \cite{HL:24}.
A \emph{(rooted) network} $N$ is a DAG for which $|R(N)|=1$, i.e., $N$
has a unique root $\rho\in V(N)$. A \emph{(rooted) tree} is a network that does not contain vertices
with $\indeg_G(v)>1$.

For every vertex $v\in V(G)$ in a DAG $G$, the set of its descendant leaves
$\CC_G(v)\coloneqq\{ x\in L(G)\mid x \preceq_G v\}$
is a \emph{cluster} of $G$. We write $\mathfrak{C}_G\coloneqq\{\CC_G(v)\mid v\in V(G)\}$ for the set
of all clusters in $G$. Note that $\mathfrak{C}_G$ is a grounded set system on
$L(G)$ for every DAG $G$ \cite{HL:24}. Moreover, $\mathfrak{C}_N$ is a clustering system for every network $N$; 
cf.\ \cite[Lemma 14]{Hellmuth2023}.

An edge $e=(u,w)$ in a DAG $G$ is a \emph{shortcut} if there is a vertex $v\in V(G)$ such that $w \prec_G v \prec_G u$. 
A DAG without shortcuts is \emph{shortcut-free}. We denote with $G^-$ the \REV{the directed graph} obtained from $G$ by removing all
shortcuts \REV{from the edge set of $G$.  The next result guarantees that removing some or all shortcuts of $G$ results in a DAG that 
preserves both the partial order $\preceq_G$ and the set system $\mathfrak{C}_G$ of the original DAG. In particular, $G^-$ is uniquely determined and shortcut-free. }
\begin{lemma}[{\cite[L.~2.5]{HL:24}}]\label{lem:shortcutfree}
	Let $G = (V,E)$ be a DAG on $X$ and $e$ be a shortcut in
	$G$. \REV{Then. $G' \coloneqq (V , E \setminus \{e\})$ is a DAG on $X$.}
	Moreover, for all $u,v\in V$, it holds
	that $u\prec_G v$ if and only if $u\prec_{G'} v$ and, for all $v\in V$, it holds that
	$\CC_G(v)= \CC_{G'}(v)$.
	In particular, $G^-$ is \REV{a shortcut-free DAG that is} uniquely determined.
\end{lemma}

For a given a DAG $G$ and a subset $A\subseteq L(G)$, a vertex $v\in V(G)$ is a \emph{common
ancestor of $A$} if $v$ is ancestor of every vertex in $A$. Moreover, $v$ is a \emph{least common
ancestor} (LCA) of $A$ if $v$ is a $\preceq_G$-minimal vertex that is an ancestor of all vertices in
$A$. The set $\LCA_G(A)$ comprises all LCAs of $A$ in $G$. 
We will, in particular, be interested in situations \REV{where} $|\LCA_G(A)|=1$ holds for certain subsets
$A\subseteq X$. For simplicity, we will write $\lca_G(A)=v$ in case that $\LCA_G(A)=\{v\}$  and say that
\emph{$\lca_G(A)$ is well-defined}; otherwise, we leave $\lca_G(A)$
\emph{undefined}. 

For later reference, we provide a connection between clusters and \REV{LCAs}.
\begin{lemma}[\REV{\cite[L.~1 \& Obs.~2]{SCHS:24}}]  
  \label{lem:deflcaY}
  Let $G$ be a DAG on $X$, $\emptyset\ne A\subseteq X$, and
  suppose $\lca_G(A)$ is well-defined. Then the following is satisfied:
  \begin{enumerate}[label=(\roman*)] 
  \item $\lca_G(A)\preceq_{G} v$ for all $v$ with $A\subseteq\CC_G(v)$.
  \item $\CC_G(\lca_G(A))$ is the unique inclusion-minimal cluster in
    $\mathfrak{C}_G$ containing $A$.
  \end{enumerate}
  Moreover, if $u\preceq_G v$, then $\CC_G(u) \subseteq \CC_G(v)$ for all $u,v\in V(G)$.
\end{lemma}

Lemma~\ref{lem:deflcaY} shows that if two vertices are $\preceq_G$-comparable, then their respective clusters are comparable with respect to inclusion. 
The following property ensures the converse, namely, $\preceq_G$-comparability
of vertices $u$ and $v$ based on subset-relations between the corresponding clusters $\CC_G(u)$
and $\CC_G(v)$. 

\begin{definition}[{\cite{Hellmuth2023}}]\label{def:PCC}
	A DAG $G$ has the path-cluster-comparability (PCC) property if it satisfies, for all $u, v \in V(G)$:
	$u$ and $v$ are $\preceq_G$-comparable if and only if $\CC_G(u) \subseteq  \CC_G(v)$ or  $\CC_G(v) \subseteq  \CC_G(u)$.
\end{definition}

As we shall see soon,
for every grounded set system $\mathfrak{C}$ there is a 
DAG $G$ with $\mathfrak{C}_G = \mathfrak{C}$ that satisfies (PCC). 
This result builds on the concepts of Hasse diagrams and regular DAGs. 	
The \emph{Hasse diagram} $\Hasse(\mathfrak{C})$ of a set system $\mathfrak{C}\subseteq 2^X$ is the
DAG with vertex set $\mathfrak{C}$ and directed edges from $A\in\mathfrak{C}$ to $B\in\mathfrak{C}$
if (i) $B\subsetneq A$ and (ii) there is no $C\in\mathfrak{C}$ with $B\subsetneq C\subsetneq A$. We
note that $\Hasse(\mathfrak{C})$ is also known as the \emph{cover digraph} of $\mathfrak{C}$
\cite{Baroni:05,sem-ste-03a}. 

\begin{definition}\label{def:regular}
A DAG $G=(V,E)$ is \emph{regular} if the map $\varphi\colon V\to
V(\Hasse[\mathfrak{C}_G])$ defined by $v\mapsto \CC_G(v)$ is an isomorphism between $G$ and
$\Hasse(\mathfrak{C}_G)$ \cite{Baroni:05}.  
\end{definition}

In general, we are interested in DAGs $G$ with certain properties and that satisfy
$\mathfrak{C}_G=\mathfrak{C}$ for a given grounded set system $\mathfrak{C}$. Structural properties
of $\Hasse(\mathfrak{C})$ are, in this context, often helpful. However,
$\mathfrak{C}_{\Hasse(\mathfrak{C})}\neq \mathfrak{C}$ holds as the leaves of $\Hasse(\mathfrak{C})$
are labeled with the inclusion-minimal elements in $\mathfrak{C}$, i.e., as sets. To circumvent
this, we write $G\doteq \Hasse(\mathfrak{C})$ for the directed graph that is obtained from
$\Hasse(\mathfrak{C})$ by relabeling all vertices $\{x\}$ in $\Hasse(\mathfrak{C})$ by $x$. Thus,
for $G\doteq \Hasse(\mathfrak{C})$ it holds that $\mathfrak{C}_G=\mathfrak{C}$ provided that
$\mathfrak{C}$ is a grounded set system on $X$.

\begin{lemma}[{\cite[Lem.~4.7]{HL:24}}]\label{lem:properties-GdotH}
	For every set system $\mathfrak{C}$, the directed graph $G\doteq \Hasse(\mathfrak{C})$ is a shortcut-free
	DAG that satisfies (PCC). 
	Moreover, if $\mathfrak{C}$ is grounded, then $\Hasse(\mathfrak{C})$
	is regular and phylogenetic. 
\end{lemma}


\paragraph{The $\boldsymbol{\ominus}$-Operator and  \Olca-\rel DAGs.}

\REV{We now focus on summarizing relevant definitions and results from \cite{HL:24} that we need here.}
\begin{definition} \label{def:LCAlcarel}
	Let $G$ be a DAG on $X$ and $\mathscr{I}$ be a set of integers. 
	A vertex	$v$ is an \emph{$\mathscr{I}$-$\lca$ vertex} if $v = \lca_G(A)$ for some $A\in X(\mathscr{I})$.
	Moreover, $G$ is \emph{\OlcaTMP-relevant} (in short \OlcarelTMP)	if  all vertices in $V(G)$ are \OlcaTMP vertices. 
\end{definition}

\REV{
\begin{observation}
	$\{1,\dots,|L(G)|\}$-\lca-\rel DAGs $G$ are precisely those DAGs where each vertex is the unique LCA for
	some subset $A\subseteq L(G)$ of leaves.
\end{observation}}

	\REV{To illustrate Definition~\ref{def:LCAlcarel},} consider the DAG $G$ in Figure~\ref{fig:exmpl-preI}(left), where 
	the vertex \REV{$v$} is not an \OlcaTMP vertex for any set of integers $\mathscr{I}$
	while the vertex \REV{$u$} is a $\{2\}$-$\lca$ and a $\{3\}$-$\lca$ vertex 
	\REV{since $u=\lca_G(A)$ for both $A=\{c,d\}$ and $A=\{b,c,d\}$. However, $u$ is} not a $\{k\}$-$\lca$ vertex 
	for any $k\neq 2,3$. Hence, the vertex \REV{$u$} is a $\mathscr{I}$-$\lca$ vertex for any set
	$\mathscr{I}$ that contains $2$ or $3$.  Due to the vertex \REV{$v$,} 
	$G$ is not \OlcarelTMP for any set  $\mathscr{I}$. 
	In contrast, the DAG in Figure~\ref{fig:exmpl-preI}(right)
	is \OlcarelTMP for any $\mathscr{I}$ with $1,2\in \mathscr{I}$.
	Note that any leaf $x$ in $G$ is a  $\{1\}$-$\lca$ but not a $\{k\}$-$\lca$ vertex 
	for any $k\neq 1$. The latter is, in particular, true for any DAG
	and, since all DAGs have at least one leaf \cite{HL:24}, there exist no DAG that is 
	$\mathscr{I}$-$\lca$-$\rel$ for any $\mathscr{I}$ with $1\notin\mathscr{I}$. 
	To address this issue, we  use

\begin{definition}\label{def:1inI}
For any DAG $G$ considered here, the set $\One$ denotes a subset of $\{1,\dots, |L(G)|\}$ that satisfies $1\in \One$. 
\end{definition}

\REV{As the next result shows, \Olcarel DAGs from which all shortcuts have been removed must be regular.}

\begin{proposition}[{\cite[Cor.~4.9]{HL:24}}]\label{prop:Onerel=>regular}
	For every \Olcarel DAG $G$, the DAG $G^-$ is regular. 
\end{proposition}

Not all DAGs are $\One$-$\lca$-\rel (cf.\ Figure~\ref{fig:exmpl-preI}). As shown in \cite{HL:24}, 
it is possible to ``transform'' a non-$\One$-$\lca$-\rel DAG $G$ into an $\One$-$\lca$-\rel
DAG $H$  while preserving essential structural
properties of the original DAG $G$ using the following $\ominus$-operator.

\begin{definition}[{\cite{HL:24}}]\label{def:ominus}
  Let $G=(V,E)$ be a DAG and $v\in V$. Then $G\ominus v=(V',E')$ is the
  directed graph with vertex set $V'=V\setminus\{v\}$ and edges $(p,q)\in E'$ 
  precisely if $v\ne p$,
  $v\ne q$ and $(p,q)\in E$, or if $(p,v)\in E$ and $(v,q)\in E$. 
  For a non-empty subset $W = \{w_1,\dots,w_\ell\} \subsetneq V$, define 
  $G\ominus W \coloneqq (\dots ((G \ominus w_1) \ominus w_2) \dots)\ominus w_\ell$.
\end{definition} 

In simple words, the directed graph $G \ominus v$ is obtained from $G$ by removing $v$ and its incident edges
and connecting each parent  $p$ of $v$ with each child $q$ of $v$. In case $v$ is a leaf or a root,  $v$ and its incident
edges are just deleted. As shown next, the operator $\ominus$ can be used to transform any DAG into an
\Olcarel version while preserving key structural properties of the original DAG.

\begin{theorem}[\cite{HL:24}]\label{thm:S4S5}
	Let $G$ be a DAG on $X$ and $W\subseteq V(G)$ be a non-empty subset of
	vertices that are not $\One$-$\lca$ vertices in $G$. 
	Then, $H\coloneqq G\ominus W$ is a DAG that satisfies 
	\begin{enumerate}[label=\emph{(S\arabic*)}]
	    \setcounter{enumi}{-1}  
 	    \item  $\mathfrak{C}_{H} \subseteq \mathfrak{C}_G$, meaning no new clusters are introduced. 
		 	    In particular,	 $\CC_G(u)=\CC_{H}(u)$ for all $u\in V(H)$.
		 \item $H$ remains a DAG on $X$, \REV{that is,} the leaves of $H$ remain the same as in $G$.
		 \item $V(H) \subseteq V(G)$, meaning no new vertices are introduced.
		 \item $H$ preserves the  ancestor relationship $\prec_G$, i.e., $u \prec_G w$ if and only if $u \prec_{H} w$ for all $u, w \in V(H)$.
		 \item $H$ is \emph{$\One$-$\lca$-preserving}, i.e., $\lca_{H}(A) = \lca_G(A)$ for all $A\in X(\One)$ for which $\lca_G(A)$ is well-defined.
	\end{enumerate}
	Moreover, if $W$ contains every vertex of $G$ that is not an \Olca vertex of $G$, then $G\ominus W$ is \Olcarel.
\end{theorem}

Hence, the $\ominus$-operator preserves essential structural features (S0) -- (S4) of $G$. 
Moreover, if $H\coloneqq G\ominus W$ is \Olcarel, then  Proposition~\ref{prop:Onerel=>regular} implies that \REV{$H^-$} is regular. 
\REV{Since, additionally, $\mathfrak{C}_H=\mathfrak{C}_{H^-}$ holds by Lemma~\ref{lem:shortcutfree}, the latter means that} $H\simeq \Hasse(\mathfrak{C}_H)$ where, by (S0), $\mathfrak{C}_{H} \subseteq \mathfrak{C}_G$ holds. 
In other words, ``removal'' of all non-$\One$-$\lca$ vertices via the $\ominus$-operator from $G$ results always in a DAG $H$ \REV{such that $H^-$} is
isomorphic to the Hasse diagram of a subset of  $\mathfrak{C}_G$. This raises the question of
how the set systems $\mathfrak{C}_G$ and $\mathfrak{C}_{G \ominus W}$ are
related.
Since $\mathfrak{C}_{G \ominus W} \subseteq \mathfrak{C}_G$ for any set $W$ consisting of non-$\One$-$\lca$ vertices in $G$,
it is  of particular interest to understand which clusters are contained
in $ \mathfrak{C}_G \setminus \mathfrak{C}_{G \ominus W}$, if there are any. As an example see Figure~\ref{fig:exmpl-preI} for the case
  $ \mathfrak{C}_G \setminus \mathfrak{C}_{G \ominus W} \neq \emptyset$ and Figure~\ref{fig:Gtree} for the case  $ \mathfrak{C}_G \setminus \mathfrak{C}_{G \ominus W} = \emptyset$.
We will answer this question for
DAGs having the so-called $\One$-lca property 
in the next section. Before, we provide the following useful lemma.

\begin{lemma}\label{lem:rel-shortcutfree}
	Let $G = (V,E)$ be a DAG on $X$,  the edge $e\in E$ be a shortcut in
	$G$ and $G' = (V,E\setminus \{e\})$.
	Then, $G'$ satisfies (S0) -- (S4).
	If $G$ is \Olcarel, then $G'$ is \Olcarel. 
\end{lemma}
\begin{proof}
	Let $G = (V,E)$ be a DAG on $X$ and $e$ be a shortcut in $G$. Put $G' = (V,E\setminus \{e\})$. 
	By Lemma~\ref{lem:shortcutfree} and since $V(G')=V$, $G'$ satisfies (S0) -- (S3) and, even more, 
	$\CC_G(v)= \CC_{G'}(v)$ for all $v\in V$. 
	To see that (S4) is satisfied,
	assume that $u\coloneqq	\lca_G(A)$ is well-defined for some  $A\in X(\One)$. 
	In particular, $u$ is the unique $\preceq_G$-minimal vertex with
	the property that $A\subseteq \CC_G(u)$. Since $G'$ is $\prec_G$-preserving and $\CC_G(w)=
	\CC_{G'}(w)$ for all $w\in V$ it follows that $u$ is the unique $\preceq_{G'}$-minimal vertex
	with the property that $A\subseteq \CC_{G'}(u)$. Hence, $u = \lca_{G'}(A)$. 
	It is now a direct consequence of (S4) that,
	if $G$ is \Olcarel, then $G'$ is \Olcarel.
\end{proof}


\section{$\OneCaption$-ary Set Systems  and  DAGs with $\OneCaption$-$\lca$-Property.}
\label{sec:ominus-one-lca-prop}

In the following, we consider a generalization of so-called lca-networks introduced in \cite{Hellmuth2023} 
and $k$-lca-DAGs introduced in \cite{SCHS:24}.

\begin{figure}
	\centering
	\includegraphics[width=0.7\textwidth]{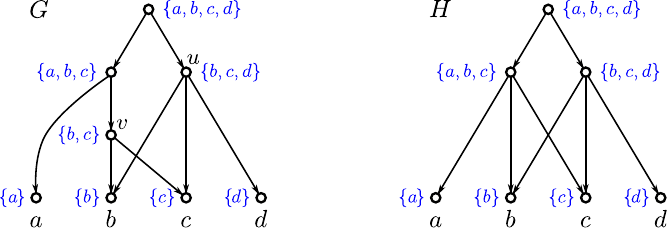}
	\caption{Shown is a network $G$ and a regular network $H = G\ominus v$. 
				\REV{Clusters $\CC(w)$ in $G$, resp., $H$ are written in blue text next to the respective vertex $w$.}
				While $H$ satisfies (PCC), $G$ does not due to the vertices \REV{$v$ and $u$  
				for which $\CC_G(v) = \{b,c\} \subseteq \CC_G(u) =\{b,c,d\}$ but neither $u\preceq_G v$ nor $v\preceq_G u$ holds.} 
			   The set $\mathfrak{C}_G$ is pre-$\{1,2\}$-ary, but $\mathfrak{C}_H$ is not	 since both cluster $\{a,b,c\}$ and $\{b,c,d\}$ are inclusion-minimal ones
			   that contain $\{b,c\}$.
			   Both, $\mathfrak{C}_G$ and $\mathfrak{C}_H$, are $\One$-ary for any $\One$
			   with $2\notin \One$. 		
			   Neither $G$ nor $H$ has $\{1,2\}$-$\lca$ property since the two leaves $b$ and $c$ do not have a unique LCA in $G$ and $H$.
			   However,	   both $G$ and $H$ have the $\One$-$\lca$-property for any $\One$
			   with $2\notin \One$. 
			   While $H$ is $\One$-$\lca$-\rel for all $\One$ containing $2$ or $3$,  
			   $G$ is not $\One$-$\lca$-\rel for any $\One$ due to vertex $v$.
	} 
	\label{fig:exmpl-preI}\label{fig:exmpl-2lca}
\end{figure}

\begin{definition}\label{def:Olca-prop}
	A DAG $G$ on $X$ has the \emph{\Olca-property} if $\lca_G(A)$ is well-defined for all $A\in X(\One)$. 
\end{definition}
\REV{The definition of the \Olca-property and that of DAGs that are \Olcarel look very similar. However, these two definitions are distinct concepts and neither of them implies the other.
By way of example, 
in Figure~\ref{fig:exmpl-2lca}(right) the network $H$ is, as previously examplified, $\{1,2\}$-\lca-\rel but since $\LCA_H(\{b,c\})$ contains the two common parents of $b$ and $c$, 
$H$ does not have the $\{1,2\}$-\lca-property. Moreover, Figure~\ref{fig:Gtree}(left) shows a DAG $G$ that can (tediously) be verified to have the $\{1,2\}$-\lca-property although it is not $\{1,2\}$-\lca-\rel. 
}

We emphasize that there is no loss of generality assuming that $1\in\One$ in the definition of the \Olca-property, since $\lca_G(\{x\})=x$ is well-defined for every leaf $x$.
A simple structural property of DAGs with \Olca-property is as follows.

\begin{lemma}\label{lem:lcaprop=>connected}
Every  DAG $G$ with the \Olca-property for some $|\One|>1$ is connected.
\end{lemma}
\begin{proof}
Let $G=(V,E)$ be a DAG on $X$ with the \Olca-property for some $\One$ such that $|\One|\geq2$.
Assume, for contradiction, that $G$ is not connected. In this case, there are two
leaves $x,y\in X$ for which there is no \REV{undirected} $xy$-path. Hence, there is no 
common ancestor of $\{x,y\}$. This, in particular, remains true for any subset $A\in X(\One)$ 
that contains $x$ and $y$ in which case $\lca_G(A)$ is not well-defined; a contradiction. 
\end{proof}

The requirement that $|\One|>1$ cannot be omitted from Lemma~\ref{lem:lcaprop=>connected}. 
To see this, consider the DAG $G = (\{x,y\}, \emptyset)$ which has the $\{1\}$-$\lca$-property
but it is not connected.
We now show that the \Olca-property is preserved under the $\ominus$-operator. 

\begin{lemma}\label{lem:Olca-property}
	Let $G$ be a DAG with the \Olca-property and $W$ \REV{be some set of} vertices that are not \Olca
	vertices of $G$. Then $G\ominus W$ has the \Olca-property.
\end{lemma}
\begin{proof}
	Since $G$ has the \Olca-property, $\lca_G(A)$ is well-defined for each $A\in X(\One)$. By
	Theorem~\ref{thm:S4S5}, $G\ominus W$ satisfies (S4) and, therefore, $\lca_{G\ominus W}(A)=\lca_G(A)$ is well-defined for each $A\in
	X(\One)$. Consequently, $G\ominus W$ has the \Olca-property.
\end{proof}

We now examine how the structural properties of a DAG $G$ that are \Olcarel or possess the
\Olca-property relate to the properties of the associated set system $\mathfrak{C}_G$.
First note that in a DAG with the 
$\One$-$\lca$-property, Lemma~\ref{lem:deflcaY}(ii) implies that $\CC_G(\lca_G(A))$ is the unique
inclusion-minimal cluster in $\mathfrak{C}_G$ containing $A$ for every $A\in X(\One)$. Hence
$\mathfrak{C}_G$ is pre-$\One$-ary.  	

\begin{observation}\label{obs:Ilca->pre-I-ary}
	For every DAG $G$ with $\One$-$\lca$-property,  the set system $\mathfrak{C}_G$ is pre-$\One$-ary. 
\end{observation}

The converse of Observation\ \ref{obs:Ilca->pre-I-ary} is, in general, not satisfied as
shown by the DAG $G$ in Figure~\ref{fig:exmpl-preI}(left).
In this example, the clustering system $\mathfrak{C}_{G}$ is pre-$\{1,2\}$-ary
but $G$ does not have the $\{1,2\}$-$\lca$-property since $|\LCA_G(\{b,c\})|>1$.
Note also that \REV{being} \Olcarel does, in general, not imply  that 
$\mathfrak{C}_G$ is pre-$\One$-ary, see Figure \ref{fig:exmpl-2lca}(right). 
In this example, the network $H$ is $\{1,2\}$-$\lca$-\rel. However,  both clusters
$\{a,b,c\}$ and $\{b,c,d\}$ in $\mathfrak{C}_{H}$ are inclusion-minimal clusters containing
$\{b,c\}$ and, thus, $\mathfrak{C}_{H} $ is not \REV{pre-$\{1,2\}$-ary.}
Nevertheless, we obtain the following result which generalizes \REV{\cite[Prop.\ 3]{SCHS:24}.
This result is of particular interest, as it shows that the property of a DAG $G$ having the \Olca-property 
is fully determined by the structure of its underlying set system $\mathfrak{C}_G$, 
provided that $G$ satisfies (PCC). }

\begin{theorem}\label{thm:PCC-Olca-preO}
Let $G$ be a DAG that satisfies (PCC). 
Then, $G$ has the \Olca-property if and only if $\mathfrak{C}_G$ is pre-$\One$-ary. 
\end{theorem}
\begin{proof}
The \emph{only if}-direction holds by Observation\ \ref{obs:Ilca->pre-I-ary}. 
Suppose now that $G$ is a DAG on $X$  that satisfies (PCC) and 
such that $\mathfrak{C}_G$ is pre-$\One$-ary. Let $A\in X(\One)$. 
We show that $\lca_G(A)$ is well-defined.
Since  $\mathfrak{C}_G$ is pre-$\One$-ary, 
there is a unique inclusion-minimal element $\CC_G(z)\in \mathfrak{C}_G$ 
such that $A\subseteq  \CC_G(z)$. Thus, $z$ is a common ancestor of $A$
which, in particular, implies that 
$\LCA_G(A)\neq \emptyset$. Assume, for contradiction, that 
$|\LCA_G(A)|\neq 1$ and let $u,w\in \LCA_G(A)$. Since $\CC_G(z)$ is 
the unique inclusion-minimal element in $\mathfrak{C}_G$ that contains $A$, 
$\CC_G(z) \subseteq \CC_G(u)$ and $\CC_G(z) \subseteq \CC_G(w)$ must hold. 
Hence, $\CC_G(z) \subseteq \CC_G(w)\cap \CC_G(u)$. 
Since both $u$ and $w$ are $\preceq_G$-minimal ancestors of the vertices in $A$, 
neither $u\prec_G w$ nor $w\prec_G u$ can hold. 
Since $G$ satisfies (PCC), neither $\CC_G(w) \subseteq \CC_G(u)$ nor $\CC_G(u) \subseteq \CC_G(w)$
can hold. This together with $\CC_G(w)\cap \CC_G(u)\neq \emptyset$ implies that  $\CC_G(w)$ and $\CC_G(u)$ must overlap
and, in particular, $\CC_G(z) \subsetneq \CC_G(u)$ and $\CC_G(z) \subsetneq \CC_G(w)$.
This and the fact that $G$ satisfies (PCC) implies that $u$ and $z$ are $\preceq_G$-comparable. 
If, however, $u\preceq_G z$, then $\CC_G(u) \subseteq \CC_G(z)$ (cf.\ Lemma~\ref{lem:deflcaY}).
Thus, $z\preceq_G u$ and, by similar arguments, $z\preceq_G w$  must hold. 
This, however, contradicts the fact that $u$ and $w$ are $\preceq_G$-minimal ancestors of the vertices in $A$. 
Consequently, $|\LCA_G(A)| = 1$ must hold and $G$ satisfies the \Olca-property.
\end{proof}

We now provide characterizations of grounded pre-$\One$-ary and $\One$-ary  set systems in terms
of the DAGs from which they derive, thereby generalizing results established in \cite{SCHS:24}.
To this end, we provide first
\begin{lemma}\label{I-rel=>I-ary}
	If $G$ is an \Olcarel DAG with the $\One$-$\lca$-property, then $\mathfrak{C}_G$ is an $\One$-ary set system.
\end{lemma}
\begin{proof}
	Let $G=(V,E)$ be an \Olcarel DAG on $X$ with $\One$-$\lca$-property. By Obs.\
	\ref{obs:Ilca->pre-I-ary}, $\mathfrak{C}_G$ is pre-$\One$-ary. To show that $\mathfrak{C}_G$ is
	$\One$-ary, let $C \in \mathfrak{C}_G$. Since $\mathfrak{C}_G$ is the set of clusters $\CC_G(w)$
	with $w\in V$, it follows that there is a vertex $w\in V$ with $C = \CC_G(w)$. Moreover, since $G$
	is \Olcarel, there is some $A\in X(\One)$ such that $w=\lca_G(A)$. By Lemma~\ref{lem:deflcaY}(i), 
		for all $v\in V$ with $A\subseteq \CC_G(v)$ it holds that $w\preceq_G v$ an
	thus, $C = \CC_G(w)\subseteq \CC_G(v)$ (cf.\ Lemma~\ref{lem:deflcaY}). Hence $C = \CC_G(w)$ is the
	unique inclusion-minimal cluster in $\mathfrak{C}_G$ that contains $A$, i.e., $C$ satisfies
	($\One$-$C$). Thus, $\mathfrak{C}_G$ is $\One$-ary.
\end{proof}

We are now in the position to provide the aforementioned characterizations. 

\begin{theorem}\label{thm:pre-I-ary_char} \label{thm:one-ary<=>OLCA}
Let  $\mathfrak{Q}$ be a set system on $X$. 
Then, $\mathfrak{Q}$ is pre-$\One$-ary and grounded if and only if there is a DAG $G$ on $X$
with \Olca-property and $\mathfrak{C}_G=\mathfrak{Q}$.
Moreover, $\mathfrak{Q}$ is $\One$-ary and grounded if and only if there is an \Olcarel
DAG $G$ on $X$ with the \Olca-property such that $\mathfrak{C}_G=\mathfrak{Q}$.
\end{theorem}
\begin{proof}
	Suppose that  $\mathfrak{Q}$ is a pre-$\One$-ary and grounded set system on $X$.
	Consider the DAG $G\doteq \hasse(\mathfrak{Q})$. 
	By definition, $G$ is a DAG on $X$ and $\mathfrak{C}_G = \mathfrak{Q}$.
	By Lemma~\ref{lem:properties-GdotH}, $G$ satisfies (PCC). 
	By Theorem \ref{thm:PCC-Olca-preO}, $G$ has the \Olca-property. 
	Conversely, suppose that $G$ is a DAG on $X$ with \Olca-property and $\mathfrak{C}_G=\mathfrak{Q}$.
	Clearly, $\mathfrak{Q}$ is thus a grounded set system on $X$.
	By Observation\ \ref{obs:Ilca->pre-I-ary}, $\mathfrak{Q}$ is pre-$\One$-ary. 

		Suppose now that $\mathfrak{Q}$ is an $\One$-ary and grounded set system on $X$.
		Consider the DAG $G\doteq \hasse(\mathfrak{Q})$. 
		By definition, $G$ is a DAG on $X$ and $\mathfrak{C}_G = \mathfrak{Q}$.
		By Lemma~\ref{lem:properties-GdotH}, $G$ satisfies (PCC).
		Since $\mathfrak{Q}$ is $\One$-ary it is, in particular, pre-$\One$-ary and Theorem~\ref{thm:PCC-Olca-preO} 
		thus implies that $G$ has the \Olca-property. It remains to show that
		$G$ is \Olcarel. To this end, let $v$ be a vertex of $G$ and put $C\coloneqq \CC_G(v)$. Since
		$C\in \mathfrak{C}_G=\mathfrak{Q}$ and $\mathfrak{Q}$ is $\One$-ary, there exist some $A\in
		X(\One)$ such that $C$ is the unique inclusion-minimal set in $\mathfrak{Q}$ containing $A$.
		Since $G$ has the \Olca-property, $\lca_G(A)$ is well-defined and it holds, by
		Lemma~\ref{lem:deflcaY}, that $C=\CC_G(\lca_G(A))$. Note that by definition of $G$, there is
		no vertex $u\neq v$ in $G$ such that $\CC_G(u)=C$. By the latter two arguments, we can thus
		conclude that $v=\lca_G(A)$ for $A\in X(\One)$, that is, $G$ is \Olcarel.		
		Conversely, suppose that $G$ is a \Olcarel DAG with the \Olca-property such that $\mathfrak{C}_G=\mathfrak{Q}$.
		By Lemma~\ref{I-rel=>I-ary}, $\mathfrak{Q}$ is $\One$-ary. 
\end{proof}

\REV{Recall that, as Theorem~\ref{thm:S4S5} shows, any DAG $G$ can be transformed into an \Olcarel version $G\ominus W$, for which $\mathfrak{C}_{G\ominus W}\subseteq\mathfrak{C}_G$ holds. It seems, in the general case, to be difficult to determine which (if any) clusters lie in $\mathfrak{C}_G\setminus\mathfrak{C}_{G\ominus W}$. In the remaining part of this section, we shall see that this set $\mathfrak{C}_{G\ominus W}$ can be fully described, whenever $G$ satisfy the \Olca-property.}
The following subset of clusters plays an essential role in  \REV{this endeavour}.
\begin{definition}
	For a given set system $\mathfrak{C}$ on $X$, let $\mathfrak{C}(\One) \REV{\coloneqq \{C\in\mathfrak{C}\mid C \text{ satisfies ($\One$-$C$)}\}}$ denote the subset of $\mathfrak{C}$ that contains each
	cluster $C\in \mathfrak{C}$ that satisfies ($\One$-$C$). 
\end{definition}

\begin{figure}
	\centering
	\includegraphics[width=\textwidth]{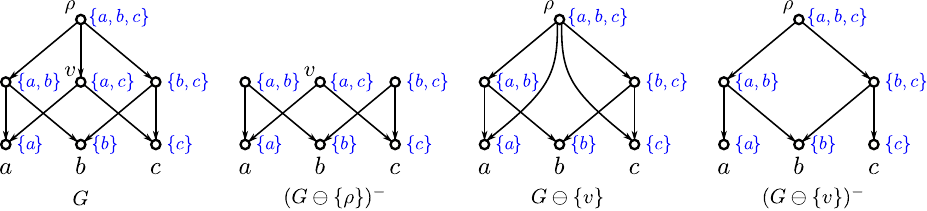}
	\caption{\REV{A DAG $G$ with the $\{1,2\}$-$\lca$-property shown alongside three DAGs obtained from $G$ by means of applying the $\ominus$-operator and removal of shortcuts. 
						The root $\rho$ of $G$ is not a $\{1,2\}$-$\lca$ vertex while all other vertices of $G$ are.
						The associated sets of clusters of these graphs examplify how $\mathfrak{C}_G(\One)$ may not be the only subset of $\mathfrak{C}_G$ that is $\One$-ary. See the text for further details.}}
	\label{fig:Exmpl-X}
\end{figure}

\REV{By definition, $\mathfrak{C}(\One)$ is a uniquely determined set system on $X$. 
Note that each cluster in  $\mathfrak{C}(\One)$ satisfies ($\One$-$C$) w.r.t.\  $\mathfrak{C}(\One)$, i.e., 
there is some $A\in X(\One)$ such that $C$ is the unique inclusion-minimal element in $\mathfrak{C}(\One)$ with $A\subseteq  C$. 
We note, however, that $\mathfrak{C}(\One)$ is not necessarily a unique subset of $\mathfrak{C}$ with this property. 
By way of example, consider the DAG $G\doteq \hasse(2^X)$ with  $X = \{a,b,c\}$ as shown in Figure~\ref{fig:Exmpl-X} and let $\One = \{1,2\}$. 
By construction, $\mathfrak{C}_G = 2^X$. It is easy to verify that}
$G$ has the \Olca-property, however, 
there is no element $A\in X(\One)=2^X\setminus \{X\}$ such that 
the cluster $X\in \mathfrak{C}_G$ is the unique inclusion-minimal cluster that contains $A$. 
Hence, $X$ does not satisfy ($\One$-$C$). In particular, 
$\mathfrak{C}_G(\One) = \mathfrak{C}_G\setminus \{X\}$ is $\One$-ary. In contrast, however, 
$\mathfrak{C}^* = \mathfrak{C}_G\setminus \{\{a,c\}\}$ is another 
subset of $\mathfrak{C}_G$ that is inclusion-maximal w.r.t. the property
of being $\One$-ary. \REV{Now consider two graphs $G\ominus \{\rho\}$ and 
$G\ominus \{v\}$ in Figure~\ref{fig:Exmpl-X}. 
We have $\mathfrak{C}_{G\ominus \{\rho\}} = \mathfrak{C}_G\setminus \{X\}$
and  $\mathfrak{C}_{G\ominus \{v\}} = \mathfrak{C}_G\setminus\{\{a,c\}\}$
and, as outlined above, both of these cluster sets are $\One$-ary. 
Hence, different vertex sets $W$ of $G$ may result in $G\ominus W$
such that $\mathfrak{C}_{G\ominus W}$ is $\One$-ary. 
In addition, both $G\ominus \{\rho\}$ and $G\ominus \{v\}$
are \Olcarel and have the \Olca-property which shows that 
the set of vertices to transform $G$ into a DAG with 
the latter two properties is, in general, not uniquely determined. }
Nevertheless, \REV{as we shall see in Theorem~\ref{thm:k-red-hasse-ominus},} the set $W$ of vertices used to transform 
a DAG $G$ with the \Olca-property into an \Olcarel DAG $G\ominus W$
that satisfies (S0) -- (S4) is uniquely determined.

We also emphasize that $\mathfrak{C}(\One)$ need not be pre-$\One$-ary and, consequently, not $\One$-ary.
For example, consider the clustering system $\mathfrak{C} = \{\{a,b,c,d\}, \{a,b,
c\},\{b,c,d\},\{a\},\{b\},\{c\},\{d\}\}$ on $X=\{a,b,c,d\}$ \REV{and put $\One\coloneqq \{1,2\}$}. Note that $\mathfrak{C} = \mathfrak{C}_H$ of the DAG $H\doteq \hasse(\mathfrak{C})$
as shown in Figure~\ref{fig:exmpl-2lca}(right). Here, \REV{$\mathfrak{C}=\mathfrak{C}(\One)$.} 
However, $\mathfrak{C}$ is not \REV{pre-$\{1,2\}$-ary}, as both
$\{a,b,c\}$ and $\{b,c,d\}$ are inclusion-minimal sets in $\mathfrak{C}$ containing $\{b,c\}\in
X(\{1,2\})$. \REV{The latter, in particular, implies that $\mathfrak{C}(\One)$ is not $\One$-ary, in general. }
As \REV{we show now in Theorem~\ref{thm:k-red-hasse-ominus},}  this situation changes \REV{for DAGs with  \Olca-property, that is, $\mathfrak{C}(\One)$ is $\One$-ary} whenever 
$\mathfrak{C}=\mathfrak{C}_G$ for some DAG $G$ with the \Olca-property.

\begin{theorem}\label{thm:k-red-hasse-ominus}
	Let $G$ be a DAG on $X$ with the \Olca-property and $W$ the set of all vertices that are not
	\Olca vertices in $G$. Then,
	 	\[(G\ominus W)^-\simeq\hasse(\mathfrak{C}_{G}(\One)).\]
	In particular, $G\ominus W$ and
	$(G\ominus W)^-$ are \Olcarel  DAGs on $X$ with the \Olca-property
	that satisfies (S0) -- (S4) and for which the set system $\mathfrak{C}_{G\ominus
	W} = \mathfrak{C}_{(G\ominus
	W)^-} = \mathfrak{C}_{G}(\One)$ is $\One$-ary. 
	\REV{In addition, $(G\ominus W)^-$  is regular.}
	Moreover, $W$ is the unique and, therefore, smallest subset of $V (G)$ such that 
	$G\ominus W$ as well as $(G\ominus W)^-$ are \Olca-\rel and
	satisfies (S0) -- (S4) w.r.t.\ $G$.
\end{theorem}
\begin{proof}
	Let $G$ be a DAG on $X$ with \Olca-property and $W$ be the set of all vertices that are not \Olca vertices in $G$. 
	Since $G$ has the \Olca-property, Observation\	\ref{obs:Ilca->pre-I-ary} implies that 
	$\mathfrak{C}_G$ is pre-$\One$-ary. 
	We continue with showing that $\mathfrak{C}_{G\ominus W} = \mathfrak{C}_G(\One)$. To this end, let first $C\in
	\mathfrak{C}_G(\One)$. Thus, $C$ satisfies ($\One$-$C$) and hence, there is some $A\in X(\One)$
	such that $C$ is the unique inclusion-minimal cluster in $\mathfrak{C}_G$ that contains $A$.
	Since $G$ has the \Olca-property, $v=\lca_G(A)$ exists and thus, $v$ is an \Olca vertex and
	therefore, not an element of $W$. By Lemma~\ref{lem:deflcaY}(ii), it holds that $\CC_G(v) = C$.
	This together with $v\notin W$ and Theorem~\ref{thm:S4S5} implies that $C
	\in\mathfrak{C}_{G\ominus W}$. Hence, $\mathfrak{C}_G(\One) \subseteq \mathfrak{C}_{G\ominus W}$.
	Now let $C\in \mathfrak{C}_{G\ominus W}$. Hence, there is some vertex $v\in V\setminus W$ in
	$G\ominus W$ such that $C=\CC_{G\ominus W}(v)$. By Theorem~\ref{thm:S4S5}, $C=\CC_{G\ominus W}(v)
	= \CC_{G}(v)$. Since $v\notin W$, $v$ is an \Olca vertex in $G$, i.e., there is some $A\in
	X(\One)$ such that $v = \lca_G(A)$. This together with Lemma~\ref{lem:deflcaY}(ii) implies that
	$C = \CC_{G}(v)$ is the unique inclusion-minimal cluster in $\mathfrak{C}_G$ that contains $A$.
	Therefore, $C\in \mathfrak{C}_G(\One)$. Hence, $\mathfrak{C}_{G\ominus W} =\mathfrak{C}_G(\One)$. 
	By Theorem~\ref{thm:S4S5}, \REV{$G\ominus W$ is \Olcarel} and, \REV{by} Lemma~\ref{lem:Olca-property}, $G\ominus W$ has the \Olca-property. 
	This allows us to apply Lemma \ref{I-rel=>I-ary} and to conclude that 
	$\mathfrak{C}_G(\One)$ is $\One$-ary. 

	\REV{We continue now with proving the remaining statements. To this end, recall first that}
	Theorem~\ref{thm:S4S5} and Lemma~\ref{lem:Olca-property} \REV{imply that} $G\ominus W$ is an \Olcarel DAG on $X$ that
	satisfies (S0) -- (S4) and that has the \Olca-property. Consider now $(G\ominus W)^-$.
	By Lemma~\ref{lem:shortcutfree}, $(G\ominus W)^-$ is uniquely determined and $\mathfrak{C}_{(G\ominus W)^-} =
	\mathfrak{C}_{G\ominus W}$. Thus, $\mathfrak{C}_{(G\ominus W)^-} $ is $\One$-ary.
	Proposition~\ref{prop:Onerel=>regular} implies \REV{that  $(G\ominus W)^-$  is regular and thus, that} $(G\ominus W)^-\simeq\hasse(\mathfrak{C}_{G\ominus W})$. 
	By \REV{the previous arguments,} $\mathfrak{C}_{G\ominus W}=\mathfrak{C}_G(\One)$ and
	it follows that 
	$(G\ominus W)^-\simeq\hasse(\mathfrak{C}_G(\One))$. Repeated application of
	Lemma~\ref{lem:rel-shortcutfree} together with the fact that $(G\ominus W)^-$ is uniquely determined
	shows that $(G\ominus W)^-$ satisfies  (S0) -- (S4).
		
	We show now that $W$ is the unique subset of $V$ such that $H \coloneqq G\ominus W$ is
	\Olcarel and satisfies (S0) -- (S4). Let $W^*$ be some subset of $V$ such that $H^*
	\coloneqq G\ominus W^*$ is \Olcarel and satisfies (S0) -- (S4). Since $G$ has the 
	\Olca-property, $\lca_G(A)$ is well-defined for all $A \in X(\One)$. 
	This together with the fact that $H$ and $H^*$ satisfy (S4) implies that 
	 \begin{equation}\label{eq:lcasame}
		\lca_G(A) = \lca_H(A) = \lca_{H^*}(A)\text{ for all }A \in X(\One).
	 \end{equation}
	In particular, this implies that $H$ and $H^*$ have the \Olca-property. Furthermore, since $H$ is
	\Olcarel, there is for each vertex $v$ of $H$ some $A\in X(\One)$ such that $v=\lca_H(A)$.
	The latter two arguments imply that $V(H)=\REV{\{\lca_H(A)\mid A\in X(\One)\}}$.
	Analogously, $V(H^*)=\REV{\{\lca_{H^*}(A)\mid A\in X(\One)\}}$. The latter two equalities together with
	Equation~\eqref{eq:lcasame} imply $V(H)=\REV{\{\lca_H(A)\mid A\in X(\One)\}}=V(H^*)$ and consequently,
	$W^* = V(G)\setminus V(H^*) = V(G)\setminus V(H) = W$. Analogous arguments show that $W$ is the
	unique subset of $V$ such that $(G\ominus W)^-$ is \Olcarel and satisfies (S0) -- (S4).
\end{proof}

We note that the converse of the first statement in Theorem~\ref{thm:k-red-hasse-ominus} is, in general, 
not satisfied, i.e.,  there are DAGs $G$ that satisfy $(G\ominus W)^-\simeq\hasse(\mathfrak{C}_{G}(\One))$
but that do not have the \Olca-property. \REV{By way of example consider the DAG $H$ in Figure~\ref{fig:exmpl-2lca}(right)
and let $\One = \{1,2\}$. In this case, $W = \emptyset$ and we obtain $ H \simeq  (H\ominus W)^-\simeq \Hasse(\mathfrak{C}_H(\One))$. 
However, $H$ does not have the \Olca-property.}
Moreover, the uniqueness of the set $W$ required to transform a DAG $G$ with the
\Olca-property into an \Olca-\rel DAG $G \ominus W$ satisfying (S0) -- (S4) is, in general, not
guaranteed for arbitrary DAGs (cf.\ e.g.\ \cite[Fig.~6]{HL:24}).

\section{DAGs with tree-like and galled-tree-like clustering systems.}
\label{sec:tree-like}

\begin{figure}
	\centering
	\includegraphics[width=0.6\textwidth]{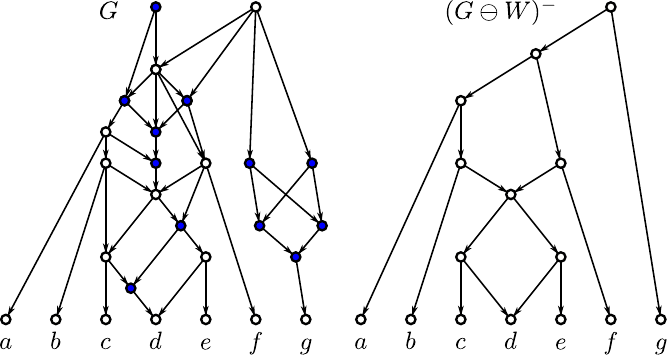}
	\caption{Shown is a DAG $G$ with $\{1,2\}$-$\lca$-property and galled-tree-like clustering system $\mathfrak{C}_G$. 
	        The vertices in the set $W$ of all non-$\{1,2\}$-$\lca$ vertices are marked in blue. According to Theorem~\ref{thm:(galled)tree-like}, 
	        the DAG $(G\ominus W)^-$ is a galled-tree that is isomorphic to $\hasse(\mathfrak{C}_G(\{1,2\}))=\hasse(\mathfrak{C}_G)$.
	} 
	\label{fig:Gtree}
\end{figure}

In what follows, we aim to determine the structure of $(G \ominus W)^-$ for DAGs $G$ 
with clustering system $\mathfrak{C}_G$ that satisfies one or more of the following properties. 
A clustering system is \emph{tree-like} if it does not contain any overlapping elements. 
Tree-like clustering systems are also known as hierarchies. 
A set system satisfies \emph{(N3O)} if it does not contain three pairwise overlapping elements
and a DAG $G$ is an \emph{(N3O)-DAG} if  $\mathfrak{C}_G$ satisfy (N3O). 
A set system $\mathfrak{C}$ satisfies property \emph{(L)} 
if $C_1\cap C_2 = C_1 \cap C_3$ for all $C_1,C_2,C_3\in \mathfrak{C}$ where $C_1$ overlaps both $C_2$ and $C_3$.
A set system $\mathfrak{C}$ is \emph{closed} if 
$A, B \in \mathfrak{C}$ and $A \cap B \neq \emptyset$  implies $A \cap B \in \mathfrak{C}$ for all $A, B \in \mathfrak{C}$. 
A network $N$ is a \emph{galled-tree} if each maximal biconnected subgraph $K$ in $N$ is either a single vertex, an edge or $K$ is composed of exactly
two directed $uv$-paths  that only have $u$ and $v$ in common \cite{Hellmuth2023}. 
A clustering system $\mathfrak{C}$ is \emph{galled-tree-like}  if $\mathfrak{C}$ is closed and satisfies (L) and (N3O). 
Note that $\mathfrak{C}_G$ can be galled-tree-like or tree-like, although $G$ is ``far away'' from being a galled-tree or tree, see 
Figure~\ref{fig:Gtree} for an illustrative example.

The following theorem, resp., lemma summarizes some results from \cite[Sec.~3.5]{sem-ste-03a} 
and \cite[Thm.~11]{Hellmuth2023}, resp., \cite[Cor.~4.11 \& L.~7.9]{HL:24}.

\begin{theorem}  \label{thm:hierarchy}
 	$\mathfrak{C}$
 	is a clustering system of a tree 
  if  and only if $\mathfrak{C}$ is tree-like.
	In this case,  $\hasse(\mathfrak{C})$
	is a  tree. 
   Moreover, $\mathfrak{C}$ is a 
  clustering system of a galled-tree if and
  only if  $\mathfrak{C}$  is galled-tree-like.
  In this case, $\hasse(\mathfrak{C})$ is a  galled-tree.
\end{theorem}

\begin{lemma}\label{lem:basics-n3o}
	For every grounded set system $\mathfrak{C}$, 
	there is a phylogenetic $\{1,\dots,|L(G)|\}$-$\lca$-\rel DAG	$G$ with $\mathfrak{C}_G = \mathfrak{C}$.
	Moreover, in an (N3O)-DAG $G$, 
	each non-leaf vertex $v\in V(G)\setminus L(G)$  is a $\{k\}$-$\lca$ vertex for some $k\geq 2$
	if and only if $v$ is an $\{\ell\}$-$\lca$ vertex for every $\ell\in\{2,\ldots,|\CC_G(v)|\}$.
\end{lemma}

As a direct consequence of Lemma~\ref{lem:basics-n3o}, we obtain 
\begin{corollary}\label{cor:(N3O)DAGexist}
	For every grounded set system $\mathfrak{C}$ that satisfies (N3O), there exists an (N3O)-DAG $G$ 
	with $\mathfrak{C}_G = \mathfrak{C}$ and such that each 
	non-leaf vertex $v\in V(G)\setminus L(G)$
	is an $\{\ell\}$-$\lca$ vertex for every $\ell\in\{2,\ldots,|\CC_G(v)|\}$.	
\end{corollary}

In the following, we call a DAG $G$ \emph{non-trivial}, if there is some  $C\in \mathfrak{C}_G$ with $|C|>1$. 
Moreover, put $\kappa_G \coloneqq \min\{|C| \colon C\in\mathfrak{C}_G,\,|C|>1\}$, i.e., $\kappa_G$ is the second smallest size of a cluster in $\mathfrak{C}_G$.

\begin{lemma}\label{lem:n3o-Cone}
	If $G$ is a non-trivial (N3O)-DAG and $\One$ is a set such that
	$k\in \One$ for some integer $k$ with $1<k\leq \kappa_G$, then $\mathfrak{C}_G = \mathfrak{C}_G(\One)$.
\end{lemma}
\begin{proof}
	Let $G$ be a non-trivial (N3O)-DAG on $X$.	Since $G$ is non-trivial, we have $\kappa_G>1$. 
	Hence, there is some integer $k$ with $1<k\leq \kappa_G$. Suppose that $\One$ contains $k\in \One$. 
	Put $\mathfrak{C} \coloneqq \mathfrak{C}_G$ and note that $\mathfrak{C}$ is grounded.
	
	By Corollary~\ref{cor:(N3O)DAGexist} there exist a DAG $H$ with $\mathfrak{C}_H = \mathfrak{C}$ and
	such that each non-leaf vertex $v$ of $H$
	is an $\{\ell\}$-$\lca$ vertex for every $\ell\in\{2,\ldots,|\CC_H(v)|\}$. 
	To show that $\mathfrak{C}(\One)=\mathfrak{C}$ holds, observe that, by construction, 
	$\mathfrak{C}(\One) \subseteq \mathfrak{C}$. Thus, it remains to show that $ \mathfrak{C} \subseteq \mathfrak{C}(\One)$. 
	To this end,
	let $C\in \mathfrak{C}$. Since $\mathfrak{C} = \mathfrak{C}_H$, there is a vertex $v\in V(H)$
	with $\CC_H(v) = C$. If $v$ is a leaf, then $|\CC_H(v)|=1$ and $1\in \One$ 
	implies that $\CC_H(v)\in\mathfrak{C}(\One)$. Suppose that $v$ is not a leaf. 
	As ensured by Corollary~\ref{cor:(N3O)DAGexist}, this vertex $v$ 
	is an $\{\ell\}$-$\lca$ vertex for every $\ell\in\{2,\ldots,|\CC_H(v)|\}$. By
	assumption, $k\in\{2,\ldots,|\CC_H(v)|\}$. 
	The latter two arguments imply that there is some $A\in X(\One)$ such that
	$v=\lca_H(A)$. By Lemma~\ref{lem:deflcaY}(ii), $\CC_H(v)$ is the unique inclusion-minimal cluster
	of $\mathfrak{C}_H$ containing $A$. By definition, $\CC_H(v)\in\mathfrak{C}(\One)$. 
	In summary, $\mathfrak{C}(\One)  = \mathfrak{C}$.
\end{proof}

\begin{theorem}\label{thm:(galled)tree-like}
	Let $G$ be a non-trivial DAG with $\One$-$\lca$-property where $\One$ contains an integer $k$ with
   with $1<k\leq \kappa_G$. Moreover, let 
	$W$ be the set of all vertices that are not \Olca vertices in $G$. 
	If $\mathfrak{C}_G$ is tree-like, then $(G\ominus W)^-$ is a phylogenetic tree $T$ that satisfies $\mathfrak{C}_T =  \mathfrak{C}_G$. 
	If $\mathfrak{C}_G$ is galled-tree-like, then $(G\ominus W)^-$ is a phylogenetic galled-tree $N$ that satisfies $\mathfrak{C}_N =  \mathfrak{C}_G$. 
\end{theorem}
\begin{proof}
	Lemma~\ref{lem:n3o-Cone} implies 
	$\mathfrak{C}_G = \mathfrak{C}_G(\One)$. 
	This together with Theorem~\ref{thm:k-red-hasse-ominus} implies that 
	$(G\ominus W)^-\simeq\hasse(\mathfrak{C}_G)$. 
	Since $\mathfrak{C}_G$ is grounded for every DAG $G$, Lemma~\ref{lem:properties-GdotH}
	implies that $(G\ominus W)^-$ is phylogenetic. 
	The latter arguments together with
	Theorem \ref{thm:hierarchy} imply that 
	$(G\ominus W)^-$ is a phylogenetic tree $T$ with clustering system $\mathfrak{C}_T =  \mathfrak{C}_G$
	 if $\mathfrak{C}_G$ is tree-like
	and that $(G\ominus W)^-$ is a phylogenetic galled-tree $N$ with clustering system $\mathfrak{C}_N =  \mathfrak{C}_G$ if $\mathfrak{C}_G$ is galled-tree-like.
\end{proof}

\REV{
\section{Summary and Outlook}

Here, we have investigated the relationship between clusters and least common ancestors (LCAs)
in directed acyclic graphs (DAGs), with a particular focus on the interplay between
$\mathscr{I}$-$\lca$-relevant DAGs and those satisfying the $\mathscr{I}$-$\lca$-property. Building
on recent work, we employed the $\ominus$-operator as a simple yet powerful tool to transform
arbitrary DAGs into $\mathscr{I}$-$\lca$-relevant ones, thereby reducing unnecessary complexity
while preserving essential structural features. For DAGs  with the $\mathscr{I}$-$\lca$-property, we
characterized precisely which clusters in $\mathfrak{C}_G$ are preserved in the simplified DAG $H =
G \ominus W$. In this setting, the vertex set $W$ to be removed is uniquely determined, ensuring a
well-defined transformation. Moreover, when $\mathfrak{C}_G$ corresponds to the clustering system of
a tree, resp., galled-tree, the transformation preserves all clusters (i.e., $\mathfrak{C}_H =
\mathfrak{C}_G$), and $H^-$ itself is always a tree, resp., galled-tree.

Several methods have been developed to simplify DAGs or networks in various ways 
\cite{Heiss2024,FRANCIS2021107215,HUBER201630,HL:24,BORDEWICH2016114}. 
One example is the normalization procedure \cite{BORDEWICH2016114,FRANCIS2021107215}, 
which is based on the notion of \emph{visible} vertices. 
A vertex $v$ in a network $N$ on $X$ is \emph{visible} if there exists a leaf $x \in X$ such that every path 
from the root of $N$ to $x$ traverses $v$.
Unfortunately, there is no direct connection between visible vertices and $\One$-$\lca$-vertices. 
For example, consider the graph $G$ shown in Figure~\ref{fig:Exmpl-X} and let $\mathscr{I} = \{1,2\}$. 
In this case, the root $\rho$ of $G$ is visible, yet it does not serve as the unique LCA for any subset $A$ 
of leaves with $|A| \in \mathscr{I}$. 
In contrast, the three other internal vertices (distinct from $\rho$) are not visible but do serve as unique LCAs 
for some subset $A$ of leaves of size $2$. 
In particular, this example also shows that even if we choose $\mathscr{I} = \{1,\dots,|X|\}$ (with $|X|=3$ in the example), 
there may still be non-visible vertices that serve as LCAs for some subset of leaves. 

However, there are cases where the normalization procedure \cite{FRANCIS2021107215} 
and our $\ominus$-operator yield identical results. 
It would therefore be interesting to investigate in more detail for which classes of DAGs 
the normalization procedure and the $\ominus$-operator produce the same simplified graph. 
Similarly, it would be worthwhile to explore under which circumstances other simplification procedures, 
such as the one proposed in \cite{Heiss2024} based on so-called LSA-vertices, 
coincide with the $\ominus$-operator.

}

\bibliographystyle{spbasic}
\bibliography{ILCA}

\end{document}